\documentclass[12pt,a4paper]{article}
\usepackage{amsmath,amsthm,amssymb}
\newcommand{\Z}{\mathbb Z}

\newcommand{\C}{\mathbb C}

\newcommand{\N}{\mathbb N}

\newcommand{\algg}{\mathfrak g}
\newcommand{\algh}{\mathfrak h}
\newcommand{\algn}{\mathfrak n}

\newcommand{\algb}{\mathfrak b}

\newcommand{\univ}{{\cal U}}
\newcommand{\vac}{{\bf 1}}
\newcommand{\sptilde}{\tilde{\ }}
\newcommand{\suchthat}{\ |\ }
\newcommand{\leadt}{\operatorname{\ell \!\text{{\it t\,}}}}
\newcommand{\NO}{\,{\raise0.25em\hbox{$\mathop{\hphantom{\cdot}}%
\limits^{_{\circ}}_{^{\circ}}$}}\,}
\newtheorem{theorem}{Theorem}[section]
\newtheorem{proposition}{Proposition}[section]

\newtheorem{lemma}{Lemma}[section]
\newtheorem{remark}{Remark}[section]
\DeclareMathOperator{\Res}{Res} \DeclareMathOperator{\End}{End}

 \DeclareMathOperator{\id}{id}
 \DeclareMathOperator{\ch}{ch}
\DeclareMathOperator{\Aut}{Aut} 
\DeclareMathOperator{\mymod}{mod}
\DeclareMathOperator{\wt}{wt}
\newcommand{\sltwofin}{{\mathfrak{sl}}(2,\C)}
\newcommand{\sltwo}{{\mathfrak{sl}}(2,\C)\sptilde}
\newcommand{\slthree}{{\mathfrak{sl}}(3,\C)\sptilde}
\title{Twisted ${\mathfrak{sl}}(3,\C)\sptilde$-modules and combinatorial identities}
\author{Ivica Siladi\'c}
\date{}
\begin{document}
\maketitle
\begin{abstract}
The main result of this paper is a combinatorial description of a basis of standard level 1 module for the twisted
affine Lie algebra $A_2^{(2)}.$ This description also gives two new combinatorial identities of G\"ollnitz (or
Rogers--Ramanujan) type. Methods used through the paper are mainly developed by J. Lepowsky, R. L. Wilson, A.
Meurman and M. Primc, and the crucial role in constructions plays a vertex operator algebra approach to standard 
representations of affine Lie algebras.
\end{abstract}
\footnotetext[2]{
Partially supported by the
Ministry of Science
and Technology of the Republic of Croatia, grant 037002.}
\footnotetext[3]{
{\em keywords:} twisted affine Lie algebras, standard modules, vertex operator algebras, 
colored partitions, Rogers-Ramanujan identities}

\section{Introduction} The results contained in this paper are another step in a try to understand the
combinatorics behind description of bases of standard modules for affine Lie algebras. The first step was made
by J. Lepowsky and R. L. Wilson in \cite{LW}, and their ideas are expanded after by A. Meurman and M.
Primc in \cite{MP3}, \cite{MP1} and \cite{MP2}.

Lepowsky and Wilson in \cite{LW} proved classical Rogers--Ramanujan identities by using representations
of affine Lie algebra $\sltwo.$ They obtained these identities by expressing the principal
character of vacuum spaces for the principal Heisenberg subalgebra of $\sltwo$ in two ways. First way is to
calculate principally specialized Weyl--Kac character formula for level 3 standard $\sltwo$-modules and this
gives product side of identities. For sum sides, they have used vertex operator construction of bases of
level 3 standard $\sltwo$-modules parametrized by partitions satisfying difference 2 conditions.
They have also constructed combinatorial bases for modules of higher levels but without proof of
linear independence, which was given after by Meurman and Primc in \cite{MP3}.

Capparelli used similar ideas in \cite{C}, this time for twisted affine Lie algebra $A_2^{(2)}$ in principal
picture and its level $3$ standard modules, and he obtained two new Rogers--Ramanujan type identities which
were not known before.

Meurman and Primc gave in \cite{MP1} a combinatorial description of bases of all standard $\sltwo$-modules
by using their ``annihilating fields''. These are obtained as vertex
operators associated to elemens from some finite--dimensional $\sltwofin$-module $R$ in the vertex
operator algebra $N(k\Lambda_0).$ For Verma module $M(\Lambda)$ with integral dominant 
$\Lambda,$ the maximal submodule $M^1(\Lambda)$ can be written as
\begin{equation}\label{intQuot}
M^1(\Lambda)=\bar{R}M(\Lambda),\ \ L(\Lambda)=M(\Lambda)/M^1(\Lambda),
\end{equation}
where $\bar{R}$ is the space of all coefficients of vertex operators associated to elements from $R.$
The idea they have used was to give a combinatorial description of the character of standard module $L(\Lambda)$
from the (easily described) character of Verma module $M(\Lambda)$ and a combinatorial basis of $M^1(\Lambda)$
obtained from (\ref{intQuot}).

If $v_\Lambda\in M(\Lambda)$ is a highest weight vector, then a spanning set for $M^1(\Lambda)$ can
be obtained from (\ref{intQuot}) in a form
\begin{equation}\label{intSpSet}
rx_1x_2\cdots x_n v_\Lambda,\ \ r\in\bar{R},\ x_i\in\sltwo.
\end{equation}
By fixing a basis for $\sltwo$ we may construct the corresponding monomial basis ${\cal P}$ of the symmetric
algebra $S(\sltwo)$. Elements from ${\cal P}$ parametrize a monomial basis in the universal enveloping
algebra $U(\sltwo)$ such that $\pi=x_1\cdots x_n\in{\cal P}$ corresponds to a basis element 
$u(\pi)=x_1\cdots x_n\in U(\sltwo).$ The elements $\pi\in{\cal P}$ are interpreted as colored partitions.

A basis for the space of ``relations'' $\bar{R}$ can also be parametrized with elements from ${\cal P}$ 
in a way that for $r\in\bar{R}$ we may write
$$r=u(\rho)+\sum_{\pi\succ\rho}c_\pi u(\pi),$$
where $\preccurlyeq$ is a suitable linear order on ${\cal P}$, and $\rho$ is the smallest partition that appears
in the corresponding vertex operator coefficient. Partition $\rho$ is called the {\em leading term} of $r$ and is
denoted by $\leadt(r)=\rho.$ A basis for $\bar{R}$ is then the set 
$\lbrace r(\rho)\suchthat \rho\in\leadt(\bar{R})\rbrace,$ where $\leadt(r(\rho))=\rho.$

The spanning set (\ref{intSpSet}) can then be reduced to the spanning set of the form
\begin{equation}\label{intRedSpan}
r(\rho)u(\pi)v_\Lambda,\ \ \rho\in\leadt(\bar{R}),\ \pi\in{\cal P}.
\end{equation}
Further reducing of this spanning set to a basis for $M^1(\Lambda)$ involves finding a linear relations between
vectors of the form $r(\rho)u(\pi)v_\Lambda,$ and these are called {\em relations among relations}. Meurman
and Primc explicity constructed required relations, thus obtaining combinatorial description of a basis
of $M^1(\Lambda).$

In \cite{MP2} they have used similar ideas to construct combinatorial basis for the basic $\slthree$-module.
One of the difficulties that they have encountered was the lack of finding enough relations among relations
for vectors (\ref{intRedSpan}) parametrized by leading terms $\leadt(\bar{R}),$ and they were forced
to use a kind of Gr\"obner base theory to correctly describe a basis for $M^1(\Lambda).$

In this work I tried to apply the same ideas to the twisted affine Lie algebra $A_2^{(2)}$ in 
homogeneous picure and its basic module. Based on the work of J. Borcea (\cite{Bo}), we get analogue of 
annihilating fields in twisted case. Borcea started from the universal vertex operator algebra $N(k\Lambda_0)$
associated to the (untwisted) affine Lie algebra $\tilde{\algg}$ (in our case $\tilde{\algg}=\slthree)$ and by 
using VOA module theory he described relations $\bar{R}^\sigma$ of standard modules for twisted Lie algebra 
$\tilde{\algg}[\sigma]$ (in our case $\tilde{\algg}[\sigma]=A_2^{(2)}$) as coefficients of twisted 
vertex operators defined on twisted module for VOA $N(k\Lambda_0)$.
By almost literally applying Meurman and Primc technique we can find a ``correct'' ordered basis of $\algg$ 
(subsection 2.2), the set of leading terms $\leadt(\bar{R}^\sigma)$ (Lemma 2.1) and a spanning set of 
the form (\ref{intRedSpan}) for the maximal submodule $M^1(\Lambda)$ of Verma module $M(\Lambda)$ 
for twisted algebra $A_2^{(2)}.$

On the other hand, relations among relations were impossible to construct using similar procedures as in
\cite{MP1} and \cite{MP2}. The problem is that the space of required relations among relations carries 
no representation structure for finite dimensional $\algg$ and therefore it is very difficult to analyze.

The solution was to generate relations among relations as coefficients of tensor product of vertex operators
associated to elements from the kernel of the map
$$\Phi:N(k\Lambda_0)\otimes \bar{R}\vac\rightarrow N(k\Lambda_0),\ \ \Phi(u\otimes v)=u_{-1}v.$$

Primc described in \cite{P} and \cite{P2} the structure of finite dimensional $\algg$-submodules of
$N(k\Lambda_0)\otimes N(k\Lambda_0)$ that generate kernel of the map $\Phi,$ and moreover he obtained
sufficient conditions for corresponding vertex operator coefficients to represent all required 
relations among relations. 

Applying this results to our special case, we encounter the problem of determining the injectivity of the map
\begin{equation}\label{intInj}
u\otimes v\mapsto \Res_z z^{n+\wt(u\otimes v)} Y^\sigma(u,z)\otimes Y^\sigma(v,z)
\end{equation}
defined on certain subspace of $\ker\Phi$, 
for an arbitrary $n\in\frac{1}{T}\Z$ (the order of automorphism is $|\sigma|=T$). By using loop module structure
of $\bar{R}^\sigma$ we are able to determine the structure of the space of coefficients (\ref{intInj}), 
which finally leads to a calculation of required dimension of relations among relations (subsection 2.3).

As a final result, we obtained combinatorial description of a basis of the basic module for twisted affine Lie algebra
$A_2^{(2)}$ (Theorem 2.1). Suitably specialized, this description gives sum side of partition identities (the product side comes
from Weyl--Kac character formula), one of which we state here (the other one is stated in Theorem 3.1):

\noindent{\bf Theorem A} {\em The number of partitions $(m_1,\ldots,m_r)$ of an integer $n$ into parts different from 2
such that difference between two consequtive parts is at least 5 (i.e. $m_{i+1}-m_i\geq 5$), and
\begin{eqnarray*}
& &m_{i+1}-m_i=5\Rightarrow m_{i+1}+m_i\not\equiv \pm 1,\pm 5,\pm 7 \mod 16\\ & &m_{i+1}-m_i=6\Rightarrow
m_{i+1}+m_i\not\equiv \pm 2,\pm 6 \mod 16\\ & &m_{i+1}-m_i=7\Rightarrow m_{i+1}+m_i\not\equiv \pm 3 \mod 16\\ &
&m_{i+1}-m_i=8\Rightarrow m_{i+1}+m_i\not\equiv \pm 4 \mod 16
\end{eqnarray*}
is equal to the number of partitions of $n$ into different odd parts.}

As mentioned above, the Rogers--Ramanujan identities appear in principal picture for $A_1^{(1)}$-modules of level 3, and the
Capparelli identities appear in principal picture for $A_2^{(2)}$-modules of level 3. On the other hand, the Capparelli
identities also appear in homogeneous picture for $A_1^{(1)}$-modules of level 1, and the identity given by Theorem A
appears in homogeneous picture for the $A_2^{(2)}$-module of level 1. So, in a sense, the above identity is
``two steps away'' from the Rogers--Ramanujan identities, and since differences between parts are large, it resembles
more to the G\"ollnitz identity (\cite{G}).

Many calculations in this work were made by the computer system GAP, which may be found at 
{\em http://www-gap.dcs.st-and.ac.uk}.

I wish to express my deep gratitude to professor Mirko Primc for his invaluable contribution to
results in this paper.

\subsection{Twisted affine Lie algebras and standard modules}
Let $\algg$ be a simple complex Lie algebra and $\sigma$ an automorphism of finite order $T.$ Denote by
$\varepsilon=\exp(\frac{2\pi\sqrt{-1}}{T}).$ The automorphism $\sigma$ acts semisimply on $\algg,$ so
we have a direct sum decomposition:
\begin{equation}{\label{alggdecomp}}
\algg=\bigoplus_{j=0}^{T-1}\algg_{[j]},\text{ where } \algg_{[j]}=\lbrace x\in\algg\suchthat\sigma
x=\varepsilon^j x\rbrace.
\end{equation}
Let $\algh$ be a Cartan subalgebra of $\algg$ and $\langle\cdot,\cdot\rangle$ a symmetric invariant bilinear
form on $\algg.$ We assume that $\langle\theta,\theta\rangle=2$ for the maximal root $\theta$ with respect to
some fixed root system basis. Set
\begin{equation}\label{affghat}
\hat{\algg}[\sigma]=\coprod_{j\in\Z}\algg_{[j\mymod T]}\otimes t^{\frac{j}{T}}+\C c.
\end{equation}
Then $\tilde{\algg}[\sigma]=\hat{\algg}[\sigma]+\C d$ is the associated twisted affine Lie algebra with the
commutator
\begin{equation}
[x(i/T),y(j/T)]=[x,y]\left((i+j)/T\right)+\frac{i}{T}\delta_{i+j,0}\langle x,y\rangle c,\ i,j\in\Z
\end{equation}
where we set $x(i/T)=x\otimes t^{i/T}$ for $x\in\algg_{[i\mymod T]}.$ Furthermore, $c$ is central
element and $[d,x(i/T)]=(i/T)x(i/T).$ We also extend the action of $\sigma$ to $\tilde{\algg}[\sigma]$ by
\begin{equation}\label{sigmaonaff}
\sigma(x(i/T))=(\sigma x)(i/T)=\varepsilon^i x(i/T),\ \sigma c=c,\ \sigma d=d
\end{equation}
(for details see \cite{K}, for example).

Let $\tilde{\algg}[\sigma]=\algn_-\oplus\algh^e\otimes\algn_+$ be a triangular decomposition of
$\tilde{\algg}[\sigma],$ and let $\Lambda\in\algh^{e^*},$ where we denoted by $\algh^{e^*}$ the dual space of
$\algh^e=\algh\rtimes \C d.$ Then the Verma module
\begin{equation}\label{VermaDef}
M(\Lambda)=U(\tilde{\algg}[\sigma])\otimes_{U(\algh^e\oplus \algn_+)}\C
\end{equation}
is a highest weight $\tilde{\algg}[\sigma]$-module with the highest weight $\Lambda.$ Here, $\C$ is viewed as
1-dimensional $(\algh^e\oplus \algn_+)$-module defined by $\algn_+\cdot 1=0$ and $h\cdot 1=\Lambda(h)$ for
$h\in\algh^e.$

Let $M^1(\Lambda)$ be the unique maximal proper $\tilde{\algg}[\sigma]$-submodule of $M(\Lambda).$ Then
\begin{equation}\label{irrdbLdef}
L(\Lambda)=M(\Lambda)/M^1(\Lambda)
\end{equation}
is an irreducible module for $\tilde{\algg}[\sigma].$ If $\Lambda$ is a dominant integral weight, then
$L(\Lambda)$ is a {\em standard} $\tilde{\algg}[\sigma]$-module.
\subsection{Loop modules}\label{loopSection}
A $\algg$-module $M$ will be called $(\algg,\sigma)$-module if, besides the action of $\algg,$ there is an
action on $M$ of a subgroup $K\subset\Aut(\algg)$ generated by $\sigma,$ such that
\begin{equation}\label{KMODdef}
\mu\pi(g)\mu^{-1}=\pi(\mu(g))
\end{equation}
for all $g\in\algg, \mu\in K$ and where $\pi:\algg\rightarrow \End M$ is representation of $\algg.$

We shall say that $(\algg,\sigma)$-module $M$ is {\em irreducible} if $M$ is the only nonzero $K$-invariant $\algg$-submodule of
$M$ (note that $M$ need not be irreducible as $\algg$-module).

Let $M$ be any finite dimensional $(\algg,\sigma)$-module. Since $\sigma^T=\id_M,$ there is a decomposition of
$M$ into eigenspaces of $\sigma$ of the form

\begin{equation}\label{MSigmaDecomp}
M=\bigoplus_{j=0}^{T-1}M_j,\ \sigma(v)=\varepsilon^j\cdot v\text{ for }v\in M_j
\end{equation}
and it follows from (\ref{KMODdef}) that (\ref{MSigmaDecomp}) and (\ref{alggdecomp}) are compatible in the
sense that $$\pi(\algg_{[i]})v_j\in M_{(i+j)\mymod T}.$$

By using the decomposition (\ref{MSigmaDecomp}) we construct the space
$$\bar{M}=\bigoplus_{j\in\Z}M_{[j\mymod T]} \otimes t^\frac{j}{T}$$ and give it a structure of
$\tilde{\algg}[\sigma]$-module by
\begin{eqnarray}
x(n+\frac{i}{T})\cdot \left(v\otimes t^{m+\frac{j}{T}}\right)&=&(x\cdot v)\otimes t^{m+n+\frac{i+j}{T}}\nonumber\\ \
d\cdot\left(v\otimes t^{m+\frac{j}{T}}\right)&=&(m+\frac{j}{T})\cdot v\otimes t^{m+\frac{j}{T}}\label{loopAbstrDef}
\end{eqnarray}
for $x\in\algg_{[i]},v\in M_j$ and $m,n\in\Z.$ From this definition we also get $\frac{1}{T}\Z$-gradation of
$\bar{M}$: $$\bar{M}=\bigoplus_{n\in\frac{1}{T}\Z}\bar{M}(n),$$ where $\bar{M}(n)=
\lbrace d\cdot\left(v\otimes t^n\right) = n\cdot v\otimes t^n\rbrace.$ We shall say that $\bar{M}$
is a {\em loop} $\tilde{\algg}[\sigma]$-module.

\begin{proposition}\label{realLoopIrrdb} Suppose that $M$ is nontrivial finite dimensional 
irreducible $(\algg,\sigma)$-module.
Then $\bar{M}$ is an irreducible loop $\tilde{\algg}[\sigma]$-module.
\end{proposition}
\begin{proof}
Suppose that $\bar{S}$ is a nontrivial $\tilde{\algg}[\sigma]$-submodule of $\bar{M}.$ Because $\bar{S}$ is
invariant under $d\in\tilde{\algg}[\sigma],$ it also acquires $\frac{1}{T}\Z$-gradation, i.e.
\begin{equation}\label{sADecomp}
\bar{S}=\bigoplus_{n\in\frac{1}{T}\Z}\bar{S}(n)=\bigoplus_{n\in\frac{1}{T}\Z}S_n\otimes t^n.
\end{equation}
Define $$S=\sum_{n\in\frac{1}{T}\Z}S_n\subset M.$$ 

We shall prove that $S_{n+\frac{i}{T}}=M_i$ for all $n\in\Z,\ i=0,\ldots,T-1.$ Let $n\in\Z$ be an
arbitrary fixed integer. First, we will prove that $$S'_n=\bigoplus_{i=0}^{T-1}S_{n+\frac{i}{T}}$$ is closed under
the action of $\algg.$ This sum is indeed direct, because $S_{n+\frac{i}{T}}$ are mutually different
$\sigma$-eigenspaces. The fact that $S'_n$ is a $(\algg,\sigma)$-module will imply that either $S'_n$ is trivial,
or $S'_n=M,$ and in later case we will have $S_{n+\frac{i}{T}}=M_i.$

Suppose that there is $n\in\Z$ such that $S'_n$ is nontrivial. Pick a nonzero $s\in S_{n+\frac{i}{T}}$ and an arbitrary $x\in\algg.$ Let
$$x=x_0+\cdots x_{T-1}$$ be decomposition of $x$ of the form (\ref{alggdecomp}). Consider the following element
of $\tilde{\algg}[\sigma]:$
\begin{eqnarray*}
x'&=&x_0(0)+\cdots+x_{T-i-1}(\frac{T-i-1}{T})+\\ & & +
x_{T-i}(-1+\frac{T-i}{T})+\cdots+x_{T-1}(-1+\frac{T-1}{T}).
\end{eqnarray*}
Then we get that 
$$x'\cdot\left(s\otimes t^{n+\frac{i}{T}}\right)\in\bigoplus_{j=0}^{T-1}S_{n+\frac{j}{T}}\otimes t^{n+\frac{j}{T}}$$ 
implying that $x \cdot s\in
S'_n.$

Second, we have to prove that $S'_n=S'_m$ for all $n,m\in\Z.$ By asumption $\bar{S}$ is nontrivial $\tilde{\algg}[\sigma]$-submodule, 
so there exists at least one $n\in\Z$ for which $S'_n$ is nontrivial. But then
\begin{eqnarray*}
x'&=&x_0(m-n)+\cdots+x_{T-i-1}(m-n+\frac{T-i-1}{T})+\\ & & +
x_{T-i}(m-n-1+\frac{T-i}{T})+\cdots+x_{T-1}(m-n-1+\frac{T-1}{T})
\end{eqnarray*}
$$x'\cdot\left(s\otimes t^{n+\frac{i}{T}}\right)\in\bigoplus_{j=0}^{T-1}S_{m+\frac{j}{T}}\otimes t^{m+\frac{j}{T}}$$ shows
that $S'_n\subset S'_m$ for any $m\in\Z.$ Since for any nonzero $s\in S'_n$ there is $x\in\algg$ such
that $x\cdot s\neq 0,$ it follows that $S'_m\neq 0.$ Exchanging the role of $m$ and $n$
we get that also $S'_m\subset S'_n,$ so $S'_n=S'_m$ for all $n,m\in\Z.$
\end{proof}
\subsection{Colored partitions}
Let $k\in\C$ and let $\univ^\sigma=U(\hat{\algg}[\sigma])/(c-k).$ Now, $\tilde{\algg}[\sigma]$-modules of level
$k$ are also $\univ^\sigma$-modules. $\univ^\sigma$ inherits the gradation by derivation $d$ from
$U(\hat{\algg}[\sigma]).$ Let $\bar{\algg}[\sigma]=\hat{\algg}[\sigma]/\C c$ and denote by ${\cal
S}^\sigma=S(\bar{\algg}[\sigma]).$ Here $S(\cdot)$ stands for a symmetric algebra. Then ${\cal S}^\sigma$ is
a commutative graded algebra.

Let $B_{[j]},\ j=0,\ldots,T-1$ be a fixed basis for $\algg_{[j]}.$ Then $B=\cup_{j=0}^{T-1}B_{[j]}$ is a basis
for $\algg.$ Let $\preccurlyeq$ be a linear order on $B.$

We form a basis $\tilde{B}$ for $\tilde{\algg}[\sigma]$ such that
$$\bar{B}=\bigcup_{j\in\Z}B_{[j\mymod T]}\otimes t^{\frac{j}{T}},\ \tilde{B}=\bar{B}\cup \lbrace c,d\rbrace$$

Extend the order $\preccurlyeq$ on $\bar{B}$ such that $$i<j\text{ implies } x(i)\prec y(j).$$ We shall
denote by $\cal P$ the basis of the symmetric algebra ${\cal S}^\sigma$ consisting of monomials in elements from
$\bar{B}.$ Each element $\pi\in{\cal P}$ is a finite product of the form $$\pi=\prod_{i=1}^{\ell}a_i(j_i),\
a_i(j_i)\in\bar{B}.$$ We say that $\pi$ is a {\em colored partition} of degree $|\pi|=\sum_{i=1}^\ell
j_i\in\frac{1}{T}\Z$ and {\em length} $\ell(\pi)=\ell$ with {\em parts} $a_i(j_i)$ of degree $j_i$ and {\em
color} $a_i.$ The set of all colored partitions of degree $n$ and length $\ell$ is denoted by ${\cal
P}^\ell(n).$

The set $\cal P$ is a monoid with the unit element $1$, and the product of monomials $\pi,\rho\in{\cal P}$ is
denoted by $\pi\rho.$

Let $\pi=a_1(i_1)\cdots a_\ell(i_s),$ $\kappa=b_1(j_1)\cdots b_\ell(j_{s'}).$
We extend the order $\preccurlyeq$ on $\cal P$ such that for $\pi\neq \kappa$ we define
$\pi\prec\kappa$ if one of the following statements hold:
\begin{description}
\item[\rm{1.}] $\ell(\pi)>\ell(\kappa)$
\item[\rm{2.}] $\ell(\pi)=\ell(\kappa),\ |\pi|<|\kappa|$
\item[\rm{3.}] $\ell(\pi)=\ell(\kappa),$ $|\pi|=|\kappa|,$ and there is $k$, $\ell(\pi)\geq k\geq 1,$
such that $i_n=j_n$ for $\ell(\pi)\geq n>k$ and $i_k<j_k.$
\item[\rm{3.}] $(i_1,\ldots,i_s)=(j_1,\ldots,j_{s'})$ and there is $k$, $\ell(\pi)\geq k\geq 1,$
such that $a_n=b_n$ for $\ell(\pi)\geq n>k$ and $a_k<b_k.$
\end{description}
Then any subset $S\subset {\cal P}$ consisting of elements of fixed degree and length less than or equal
to some fixed $\ell$ has minimal element and $\mu\preccurlyeq\nu$ implies $\pi\mu\preccurlyeq\pi\nu.$
(cf. Lemmas 6.2.1 and 6.2.2 in \cite{MP1}).
\subsection{Leading terms}
For $n\in\frac{1}{T}\Z$ denote by ${\cal U}^\sigma(n)$ homogeneous component of ${\cal U}^\sigma$ (homogeneous
according to the gradation by $d$). Set $$V_p(n)=\sum_{i\geq p}{\cal U}^\sigma(n-i){\cal U}^\sigma(i),\
p\in\N.$$ Then $\{V_p(n),\ p\in\N\}$ forms a fundamental system of neighborhoods of $0\in{\cal U}^\sigma(n).$
Since $$\bigcap_{p=1}^\infty V_p(n)=\{0\}$$ it follows that $({\cal U}^\sigma(n),+)$ is a Hausdorff topological
group. Denote by $\overline{{\cal U}^\sigma(n)}$ its completion. Then
$$\overline{{\cal U}^\sigma}=\coprod_{n\in\frac{1}{T}\Z}\overline{{\cal U}^\sigma(n)}$$
is a topological ring (cf. \cite{FZ}).

Fix monomial basis
$$u(\pi)=b_1(j_1)\cdots b_n(j_n),\ \ \pi\in{\cal P},$$
of algebra ${\cal U}^\sigma.$ For $\pi\in{\cal P},\ |\pi|=n,$ set
$$U^{\cal P}_{[\pi]}=\overline{\C\text{-span}\{u(\pi')\suchthat |\pi'|=|\pi|,\pi'\succcurlyeq\pi\}}\subset
\overline{{\cal U}^\sigma(n)}$$
$$U^{\cal P}_{(\pi)}=\overline{\C\text{-span}\{u(\pi')\suchthat |\pi'|=|\pi|,\pi'\succ\pi\}}\subset
\overline{{\cal U}^\sigma(n)},$$ where the closure is taken in $\overline{{\cal U}^\sigma(n)}.$ Denote also
$$U^{\cal P}(n)=\bigcup_{\pi\in{\cal P},|\pi|=n}U^{\cal P}_{[\pi]},\ \
U^{\cal P}=\coprod_{n\in\frac{1}{T}\Z}U^{\cal P}(n)\subset\overline{{\cal U}^\sigma}$$ For $u\in U^{\cal
P}_{[\pi]},\ u\not\in U^{\cal P}_{(\pi)}$ we define the {\em leading term}
$$\leadt(u)=\pi.$$
Quote the following propositions from \cite{P}:
\begin{proposition}\label{leadExist}Every element $u\in U^{\cal P}(n),\ u\neq 0$ has a unique
leading term $\leadt(u).$
\end{proposition}
\begin{proposition}For all $u,v\in U^{\cal P}\setminus \{0\}$ we have
$\leadt(uv)=\leadt(u)\leadt(v).$
\end{proposition}
\begin{proposition}\label{leadSubspace}Let $W\subset U^{\cal P}$ be a
finite dimensional subspace and let $w:\leadt(W)\rightarrow W$ be a map such that
$$\rho\mapsto w(\rho),\ \ \leadt(w(\rho))=\rho.$$
Then $\{w(\rho)\suchthat \rho\in\leadt(W)\}$ is a basis of $W.$
\end{proposition}
\subsection{Leading terms of coefficients of vertex operators}
Let $k\in\C.$ For any untwisted affine Lie algebra $\tilde{\algg},$ its generalized Verma module
$$N(k\Lambda_0)=U(\tilde{\algg})\otimes_{U(\tilde{\algg}\geq 0+\C c)}\C v_k,$$
carries the structure of vertex operator algebra if $k\neq -g^\vee$ ($g^\vee$ is the dual Coxeter number of
$\algg$). This is so-called {\em vacuum} representation of $\tilde{\algg}$ at level $k.$ Vertex operator
$Y(v,z)$ associated to an element $v\in N(k\Lambda_0)$ is
$$Y(v,z)=\sum_{n\in\Z}v_n z^{-n-1},\ v_n\in\End(N(k\Lambda_0))$$
i.e. it is a formal Laurent series, and coefficients are linear operators on $N(k\Lambda_0)$ (for details see, for example,
\cite{FZ}, \cite{Li}, \cite{MP1}).

Let $\sigma$ be an automorphism of $\algg$ of order $T$ so that we have the associated twisted affine Lie
algebra $\tilde{\algg}[\sigma].$ Let $M$ be any restricted $\hat{\algg}[\sigma]$ module of level $k.$ Then
there is a map
$$Y_M^\sigma(v,z)=\sum_{n\in\frac{1}{T}\Z}v_n z^{-n-1},\ v\in N(k\Lambda_0),\ v_n\in\End(M)$$
that makes $(M,Y_M^\sigma)$ a (weak) twisted $N(k\Lambda_0)$-module (\cite{Li2}).

Due to the results in \cite{FZ} it is possible to consider coefficients of $Y(v,z)$ as
members of the completion $\overline{\cal U}$ of the universal enveloping algebra of $\tilde{\algg}.$

Likewise, the twisted vertex operators
$$x(z)=\sum_{n\in\Z}x_{n+\frac{s}{T}}z^{-n-\frac{s}{T}-1},\ x\in\algg_{[s]}$$
and their $m$-products generate a vertex operator algebra $V$ (see \cite{Li2}.) From \cite{Li2} and \cite{P} we
get that there is a surjective homomorphism $N(k\Lambda_0)\rightarrow V$ of vertex operator algebras given by
the map
\begin{equation}\label{YsigmaDef}
v\mapsto Y^\sigma(v,z)=\sum_{n\in\frac{1}{T}\Z}v_n z^{-n-1}\in \overline{{\cal
U}^\sigma}[[z^\frac{1}{T},z^{-\frac{1}{T}}]].
\end{equation}
The action of $Y^\sigma(v,z)$ on any restricted
$\tilde{\algg}[\sigma]$-module of level $k$ gives rise to the twisted representation $Y_M^\sigma(v,z)$ of
$N(k\Lambda_0)$ on $M.$

Set
$$U_{\text{loc}}^\sigma=\C\text{-span}\{v_n\suchthat v\in N(k\Lambda_0),n\in\frac{1}{T}\Z\}\subset
\overline{{\cal U}^\sigma}$$ where $v_n$ denotes a coefficient in $Y^\sigma(v,z).$ Then $U^\sigma_{\text{loc}}$
is a Lie algebra, and denote by $U^\sigma$ the associative subalgebra of $\overline{{\cal U}^\sigma}$ generated
by $U^\sigma_{\text{loc}}.$ We have:
\begin{proposition}(\cite{P}, proposition 2.6)
$${\cal U}^\sigma\subset U^\sigma\subset U^{\cal P}\subset\overline{{\cal U}^\sigma}.$$
\end{proposition}
Effectively, this means that we can calculate the leading terms of coefficients of both untwisted and twisted
vertex operators (cf. Proposition \ref{leadExist}).
\subsection{Relations on standard modules}\label{defR}
Let $k\neq -g^\vee$ be a positive integer. Let $R\in N(k\Lambda_0)$ be a $\algg$-module generated by the
singular vector $x_\theta(-1)^{k+1}\vac.$ Then $R$ is finite dimensional and invariant under the action of
$d,\sigma$ and $\coprod_{n\in\Z_{\geq 0}}\algg(n).$ 

Denote by
$$\bar{R}^\sigma=\C\text{-span}\{r_n\suchthat r\in R,n\in\frac{1}{T}\Z\}\subset U^\sigma.$$
Then $\bar{R}^\sigma$ is a loop $\tilde{\algg}[\sigma]$-module under the adjoint action
$$[x_m,r_n]=\sum_{i\geq 0}\binom{m}{i}(x_ir)_{m+n-i},\ r\in R,\ x\in\algg.$$
We shall also write $\bar{R}$ in the case $\sigma=\id.$

\begin{proposition}\label{maxByR}(\cite{Bo}, Theorem 4.11) Let $M(\Lambda)$ be a Verma module for
$\tilde{\algg}[\sigma]$ such that $\Lambda$ is dominant integral of level $k.$ Then
$$\bar{R}^\sigma M(\Lambda)=M^1(\Lambda),$$
where $M^1(\Lambda)$ is the maximal submodule of $M(\Lambda).$
\end{proposition}

By our construction, any nonzero element in the highest weight module $M(\Lambda)/\bar{R}^\sigma M(\Lambda)$ is
annihilated by all vertex operators $Y^\sigma(v,z)$ for $v\in\bar{R} N(k\Lambda_0).$ Therefore, they are called
{\em annihilating fields} and corresponding coefficients are called {\em relations} for
$M(\Lambda)/M^1(\Lambda).$

\subsection{Combinatorial description of spanning set for standard modules}
The subspace $R$ from the previous section is finite dimensional, so the loop $\tilde{\algg}[\sigma]$-module
$\bar{R}^\sigma\subset U^\sigma$ is a direct sum of finite dimensional homogeneous subspaces. By Proposition
\ref{leadSubspace}, a basis of $\bar{R}^\sigma$ can be parametrized by a set of leading terms 
$\leadt(\bar{R}^\sigma).$ Explicitly, there exists a map
$$\leadt(\bar{R}^\sigma)\rightarrow \bar{R}^\sigma,\ \rho\mapsto r(\rho)$$
such that $r(\rho)\in \bar{R}^\sigma,\ \leadt(r(\rho))=\rho$ and the set $\{r(\rho)\suchthat
\rho\in\leadt(\bar{R}^\sigma)\}$ is a basis of $\bar{R}^\sigma.$ Furthermore, we shall fix the map $r$ such
that the coefficient $C$ in the ``expansion'' $r(\rho)=C u(\rho)+\ldots$ is $C=1.$

A basis for Verma module $M(\Lambda)$ can be obtained by using Poincar\'e--Birkhoff--Witt theorem. 
Proposition \ref{maxByR} together with our parametrization of basis vectors of $\bar{R}^\sigma$ imply that 
maximal submodule $M^1(\Lambda)$ is spanned by vectors of the form
\begin{equation}\label{spanSet}
r(\rho)u(\pi)v_\Lambda,\ \ \rho\in\leadt(\bar{R}^\sigma),\pi\in{\cal P}.
\end{equation}
This spanning set is {\em not} a basis of $M^1(\Lambda)$, but it can be reduced to a basis by using
the so-called {\em relations among relations}, which are linear equations that describe linear dependence
between operators of the form $r(\rho)u(\pi).$

Since $r(\rho)=0$ on $L(\Lambda)$ it follows that 
$$u(\rho) = \sum_{\pi'\succ \rho}c_{\pi'}u(\pi')\text{ on } L(\Lambda)$$ 
allowing one to express $u(\rho)v_\Lambda$ in terms of $u(\pi')v_\Lambda,\ \pi'\succ\rho.$ By using induction
and (\ref{spanSet}) the basis for $L(\Lambda)$ can then be described as
$$u(\pi)v_\Lambda,\ \pi\in{\cal P},\ \rho\not\subset\pi,\ \rho\in\leadt(\bar{R}^\sigma)$$
where $\rho\not\subset\pi$ means that partition $\rho$ is not contained as subpartition in $\pi.$

\subsection{Linear independence theorem}
For colored partitions $\kappa,\rho$ and $\pi=\kappa\rho$ we shall write $\kappa=\pi/\rho$ and $\rho\subset\pi$
(meaning that $\rho$ is {\em contained} in $\pi$), and we shall say that $\rho\subset\pi$ is an {\em
embedding}. For an embedding $\rho\in\pi,$ where $\rho\in\leadt(\bar{R}^\sigma)$ we define the element
$u(\rho\subset\pi)\in U^\sigma$ by
$$u(\rho\subset\pi)=\left\{\begin{array}{lr}
u(\pi/\rho)r(\rho)&\text{ if } |\rho|>|\pi/\rho|\\ r(\rho)u(\pi/\rho)&\text{ if }|\rho|\leq|\pi/\rho|
\end{array}\right.$$
We shall consider a slightly modified (cf. \ref{spanSet}) spanning set of $M^1(\Lambda):$
\begin{equation}\label{spModif}
u(\rho\subset\pi)v_\Lambda,\ \ \rho\in\leadt(\bar{R}^\sigma),\ \pi\in{\cal P}.
\end{equation}
For a colored partition $\pi$ set
$$N(\pi)=\max\{{\cal E}(\pi)-1,0\},\ \ {\cal
E}(\pi)=|\{\rho\in\leadt(\bar{R}^\sigma)\suchthat\rho\subset\pi\}|$$ Define
\begin{equation}\label{defPhi}
\Phi:N(k\Lambda_0)\otimes N(k\Lambda_0)\rightarrow N(k\Lambda_0),\ \ \Phi(a\otimes b)=a_{-1}b.
\end{equation}
Note that $\Phi$ intertwines the action of $L_{-1},L_0,\algg(0)$ and $\sigma$ so $\ker\Phi$ is invariant under
the action of these operators. Set
$$(U^\sigma\bar{\otimes}U^\sigma)(n)=\prod_{i+j=n}(U^\sigma(i)\otimes U^\sigma(j)),\ \
U^\sigma\bar{\otimes}U^\sigma=\coprod_{n\in\frac{1}{T}\Z}(U^\sigma\bar{\otimes}U^\sigma)(n)$$ Let $a, b\in
N(k\Lambda_0)$ be homogeneous elements and
$$Y(a,z)=\sum_{n\in\frac{1}{T}\Z}a(n)z^{-n-\wt(a)},\ \ 
Y(b,z)=\sum_{n\in\frac{1}{T}\Z}b(n)z^{-n-\wt(a)}$$
corresponding vertex operators.
For a fixed $n\in\frac{1}{T}\Z$ define
$$(a\otimes b)(n)=\sum_{p+r=n}a(p)\otimes b(r)\in(U^\sigma\bar{\otimes}U^\sigma)(n).$$
For a subspace $Q\subset N(k\Lambda_0)\otimes N(k\Lambda_0)$ we shall also write
\begin{equation}\label{tensorCoef}
Q(n)=\left\{q(n)\suchthat q\in Q\right\}\subset (U^\sigma\bar{\otimes}U^\sigma)(n)
\end{equation}
Let $q=\sum a\otimes b$ be a homogeneous element in $N(k\Lambda_0)\otimes\bar{R}\vac.$ Then each element of the
sequence $q(n)=(\sum a\otimes b)(n)$ can be written uniquely as a finite sum of the form
$$c_i=\sum_{\rho\in\leadt(\bar{R}^\sigma)}b_\rho\otimes r(\rho).$$
Denote by $S$ the set of all partitions $\pi_i$ that are smallest (according to the $\preccurlyeq$ order)
partitions of the form $\rho\leadt(b_\rho)$ that appear in expression for $c_i.$ By assumption on the order
$\preccurlyeq,$ $S$ has a minimal element, and we will call it the leading term $\leadt(q(n))$ of $q(n).$

For a subspace $Q\subset (U^\sigma\bar{\otimes}\bar{R}^\sigma)$ and $n\in\frac{1}{T}\Z$ we shall also set
$$\leadt(Q(n))=\left\{\leadt(q(n))\suchthat q\in Q\right\}.$$

As before, for details see \cite{P}.

Finally, we shall use the observation that $N(k\Lambda_0)\otimes N(k\Lambda_0)$ posses a natural filtration
$(N(k\Lambda_0)\otimes N(k\Lambda_0))_\ell,\ \ell\in\Z_{\geq 0}$ inherited from filtration ${\cal U}_\ell$ of
universal enveloping algebra of $\hat{\algg}[\sigma].$
\begin{theorem}\label{theThm} (\cite{P}, Theorem 2.12) Let $Q\subset\ker(\Phi|(N(k\Lambda_0)\otimes\bar{R}\vac)_\ell)$
be a finite dimensional subspace and $n\in\frac{1}{T}\Z.$ Assume that $\ell(\pi)=\ell$ for all
$\pi\in\leadt(Q(n)).$ If
\begin{equation}\label{combRel}
\sum_{\pi\in{\cal P}^\ell(n)} N(\pi)=\dim Q(n)
\end{equation}
then for any two embeddings $\rho_1\subset\pi$ and $\rho_2\subset\pi$ in $\pi\in{\cal P}^\ell(n),$ where
$\rho_1,\rho_2\in\leadt(\bar{R}^\sigma),$ we have a relation
$$u(\rho_1\subset\pi)\in u(\rho_2\subset\pi)+\C\text{{\rm -span}}\overline{\{u(\rho\subset \pi')\mid
 \rho \in\leadt(\bar{R}^\sigma), \rho \subset \pi', \pi \prec \pi'\}}$$
\end{theorem}
Roughly speaking, this theorem states the following: if the total number of different embeddings of leading
terms from the set $\leadt(\bar{R}^\sigma)$ into the colored partitions of a fixed length $\ell$ equals to a
dimension of subspace of operators $Q(n)$ generated from vectors in $\ker\Phi,$ then for any two different
embeddings $\rho_1,\rho_2\subset\pi$ there exists a {\em relation among relation} between $u(\rho_1\subset\pi)$
and $u(\rho_2\subset\pi).$ The existence of relation among relations guarantees that the spanning set
(\ref{spModif}) would be a basis (up to the so-called {\em initial conditions}) if, for a fixed $\pi$, we take
only one of all $u(\rho\subset\pi),\ \rho\in\leadt(\bar{R}^\sigma).$

\section{Twisted affine Lie algebra $A_2^{(2)}$ and its standard module of level 1}
In this section we shall apply the theory developed in the previous sections to the special case of the basic
module for the twisted affine Lie algebra $A_2^{(2)}.$
\subsection{Construction of $A_2^{(2)}$}
Let $\algg$ be the simple Lie algebra $\mathfrak{sl}(3,\C).$ This is a rank 2 Lie
algebra of dimension 8. Let $e_1,e_2,f_1,f_2$ be its system of canonical generators, and denote the associated
simple roots by $\alpha_1,\alpha_2.$

There is only one nontrivial automorphism of Dynkin diagram of $\algg$: it permutes the simple roots $\alpha_1$
and $\alpha_2.$ Clearly, its order is 2. Denote by $\sigma$ the automorphism of Lie algebra $\algg$ induced by
this Dynkin diagram automorphism, defined on the generators as follows:
$$\sigma(e_1)=e_2,\ \sigma(e_2)=e_1,\ \sigma(f_1)=f_2,\ \sigma(f_2)=f_1.$$
Let
$$\algg=\algg_{[0]}\oplus\algg_{[1]}$$
be the eigenspace decomposition of $\algg$ under the action of $\sigma$: $\sigma$ acts as identity on
$\algg_{[0]}$ and as $-1$ on $\algg_{[1]}.$

Let $\algb$ be the Cartan subalgebra of $\algg_{[0]}.$ Let $\langle\cdot,\cdot\rangle$ be a nondegenerate
symmetric $\algg$-invariant bilinear form on $\algg.$ As a multiple of the Killing form,
$\langle\cdot,\cdot\rangle$ is also $\sigma$ invariant and remains nonsingular on $\algb,$ so $\algb$ and
$\algb^*$ may be identified via this form.

Let $\lbrace e_1,e_2,h_1,h_2,f_1,f_2,[e_1,e_2],[f_2,f_1]\rbrace$ be the basis of $\algg$ built from the set of
its canonical generators. Introduce the following notation (cf. \cite{K}, \cite{Bo}):
\begin{eqnarray*}
& &E_0=[f_2,f_1],\ E_1=\sqrt{2}(e_1+e_2)\ F_0=[e_1,e_2],\ F_1=\sqrt{2}(f_1+f_2)\\ & &H_0=[E_0,F_0],\
H_1=[E_1,F_1],\\ & &E_2=[E_1,E_0],\ F_2=[F_1,F_0],\ H_2=h_1-h_2\\ & &\theta=\alpha_1+\alpha_2\\ &
&\beta_0=(-\alpha_1-\alpha_2)|_\algb,\ \beta_1=\alpha_1|_\algb=\alpha_2|_\algb
\end{eqnarray*}
The sets $\lbrace E_1,H_1,F_1\rbrace$ and $\lbrace E_0,E_2,H_2,F_2,F_0\rbrace$ are linear bases for the spaces
$\algg_{[0]}$ and $\algg_{[1]}$ respectively, and they are normalized so that $[E_i,F_i]=H_i,\ i=0,1,2$ hold.
Furthermore, the space $\algg_{[0]}$ is Lie algebra ${\mathfrak {sl}}(2,\C)$ and $\algg_{[1]}$ is its
irreducible module.

For $i,j=0,1$ define $a_{ij}=\beta(H_i)$ and denote by $A$ the matrix $(a_{ij})_{i,j=0,1}.$ Then
$$A=\begin{pmatrix}2&-1\\-4&2\end{pmatrix}.$$

Matrix $A$ is a generalized Cartan matrix of affine type associated to the twisted affine Lie algebra
$A_2^{(2)}.$

Denote now by $\C[t,t^{-1}]$ the $\C$-algebra of Laurent polynomials in indeterminate $t$. Define the Lie
algebras
$$\hat{\algg}[\sigma]=\algg_{[0]}\otimes\C[t,t^{-1}]\oplus
\algg_{[1]}\otimes\C[t^\frac{1}{2},t^{-\frac{1}{2}}]\oplus\C c,\ \
\tilde{\algg}[\sigma]=\hat{\algg}[\sigma]\oplus\C d$$ by the conditions
\begin{eqnarray*}
& &[a\otimes t^{m/2},b\otimes t^{n/2}]=[a,b]\otimes t^{(m+n)/2}+\frac{m}{2}\delta_{m+n,0}\langle a,b\rangle c,\\ & &c\neq
0 \text{ central, }[d,a\otimes t^{m/2}]=\frac{m}{2}\cdot a\otimes t^{m/2},
\end{eqnarray*}
for $m,n\in\Z$ and $a\in\algg_{[m\mod 2]}, b\in\algg_{[n\mod 2]}.$ The Lie algebra $\tilde{\algg}[\sigma]$ is
{\em homogenous} realization of twisted Lie algebra $A_2^{(2)}.$

Next, we shall describe generators and root system of $\tilde{\algg}[\sigma].$ First, the space
$\algh=\algb\oplus\C c$ is a Cartan subalgebra of $\hat{\algg}[\sigma],$ and $\algh^e=\algh\oplus\C d$ is a
Cartan subalgebra of $\tilde{\algg}[\sigma].$ Let $\delta\in\algh^{e^*}$ be such that $\delta|_\algh=0,$
$\delta(d)=1.$ Define the linear functionals $\alpha_0,\alpha_1\in\algh^{e^*}$ by $\alpha_1|_\algb=\beta_1,\
\alpha_1(c)=0,\ \alpha_1(d)=0$ and $\alpha_0=\frac{\delta}{2}-2\alpha_1.$

Define the following elements of $\hat{\algg}[\sigma]$:
\begin{eqnarray*}
& &e_0=E_0\otimes t^\frac{1}{2},\ e_1=E_1\otimes t^{0}\\ & &f_0=F_0\otimes t^{-\frac{1}{2}},\ f_1=F_1\otimes
t^{0}\\ & &h_0=H_0+\frac{1}{2}c,\ h_1=H_1
\end{eqnarray*}
The set $\lbrace e_0,e_1,f_0,f_1,h_0,h_1,d\rbrace$ is a system of canonical generators of
$\tilde{\algg}[\sigma].$ Root system (and its geometry) may be visualized as follows:
\begin{center}
\special{em:linewidth 0.4pt} \unitlength 1mm \linethickness{0.4pt}
\begin{picture}(60.00,35.00)
\put(20.00,20.00){\circle*{1.33}} \put(30.00,20.00){\circle*{1.33}} \put(10.00,20.00){\circle*{1.33}}
\put(20.00,30.00){\circle*{1.49}} \put(30.00,30.00){\circle*{1.33}} \put(40.00,30.00){\circle*{1.49}}
\put(10.00,30.00){\circle*{1.49}} \put(0.00,30.00){\circle*{1.33}} \put(20.00,10.00){\circle*{1.49}}
\put(30.00,10.00){\circle*{1.49}} \put(40.00,10.00){\circle*{1.49}} \put(10.00,10.00){\circle*{1.33}}
\put(0.00,10.00){\circle*{1.33}} \put(-3.00,33.90){\makebox(0,0)[lt]{$\alpha_0$}}
\put(30.67,23.33){\makebox(0,0)[lt]{$\alpha_1$}}
\end{picture}
\end{center}
\subsection{Leading terms}

In order to describe our main result, we will introduce another set of symbols: take $X_1=F_0,\ X_2=F_2,\
X_3=H_2,\ X_4=E_2,\ X_5=E_0,\ X_6=E_1,\ X_7=H_1,\ X_8=F_1$. These vectors form a linear basis for Lie algebra
$\algg={\mathfrak {sl}}(3,\C).$ We will work simultaneously with untwisted ($\tilde{\algg}$) and twisted
($\tilde{\algg}[\sigma]$) affine Lie algebras, and it is worthwile to emphasize that we will use the same
notation for operators $X_i(k),k\in\frac{1}{2}\Z$; first time as operators on $N(\Lambda_0),$ and
second time as operators on $M(\Lambda).$ Hopefuly there will be no confusion for reader, since precise meaning
should be clear from the contex.

For example, some elements from $\tilde{\algg}[\sigma]$ (and therefore also the operators on $M(\Lambda)$) are
denoted as follows:
\begin{eqnarray*}
& &e_0=X_5\otimes t^{1/2}=X_5(\frac{1}{2}),\ e_1=X_6\otimes t^0=X_6(0),\\ & &[e_1,e_0]=[E_1,E_0]\otimes
t^{1/2}=E_2\otimes t^{1/2}=X_4\otimes t^{1/2}=X_4(\frac{1}{2})\\ & &f_0=X_1\otimes t^{-1/2}=X_1(-\frac{1}{2}),\
f_1=X_8\otimes t^0=X_8(0)\\ & &[f_1,f_0]=[F_1,F_0]\otimes t^{-1/2}=F_2\otimes t^{-1/2}=X_2\otimes
t^{-1/2}=X_2(-\frac{1}{2})\\ & &h_1\otimes t^0=X_7\otimes t^0=X_7(0), H_2\otimes t^{1/2}=X_3\otimes
t^{1/2}=X_3(\frac{1}{2})
\end{eqnarray*}
Define a linear order on the set $X_i,i=1,\ldots,8$ by
$$X_1\geq X_6\geq X_2\geq X_7\geq X_3\geq X_8\geq X_4\geq X_5.$$
The $\algg$ invariant subspace $R$ of $N(\Lambda_0)$ generated by the vector $X_1(-1)^2\vac(=F_0(-1)^2\vac)$ is
27 dimensional and has a basis consisting of the vectors (cf. section \ref{defR}):
\begin{center}
\begin{tabular}{cc}
$X_1X_6$                        &$X_1X_1                         $\\ $X_1X_7-X_2X_6$                 &$X_1X_2
$\\ $X_1X_8+X_2X_7-X_3X_6$          &$2X_1X_3+X_2X_2                 $\\ $X_2X_8+X_3X_7-X_4X_6$
&$X_1X_4+X_2X_3-X_6              $\\ $X_3X_8+X_4X_7-X_5X_6$          &$4X_1X_5+4X_2X_4+2X_3X_3+X_7    $\\
$X_4X_8+X_5X_7$                 &$X_2X_5+X_3X_4+X_8              $\\ $X_5X_8$
&$2X_3X_5+X_4X_4                 $\\
                                &$X_4X_5                         $\\
$2X_1X_7+X_2X_6+6X_1$           &$X_5X_5                         $\\ $4X_1X_8+2X_3X_6+X_2X_7+6X_2$\ \ \ \
&\\ $X_2X_8+X_4X_6+2X_3$            &\\ $2X_3X_8-X_4X_7+4X_5X_6+6X_4$   &\\ $X_4X_8-2X_5X_7+6X_5$           &\\
&\\ \multicolumn{2}{c}{$7X_6X_6-8X_1X_3+3X_2X_2$}\\ \multicolumn{2}{c}{$7X_6X_7+12X_1X_4-2X_2X_3+16X_6$}\\
\multicolumn{2}{c}{$7X_7X_7-48X_1X_5-6X_2X_4+4X_3X_3-14X_6X_8+16X_7$}\\
\multicolumn{2}{c}{$7X_7X_8-12X_2X_5+2X_3X_4+16X_8$}\\ \multicolumn{2}{c}{$7X_8X_8-8X_3X_5+3X_4X_4$}\\ &\\
\multicolumn{2}{c}{$5X_7X_7-48X_1X_5+12X_2X_4-4X_3X_3+20X_6X_8+8X_7$}\\
\end{tabular}
\end{center}
where $X_1X_6$ actualy stands for $X_1(-1)X_6(-1)\vac$ and so forth. Note that this space decomposes as a
direct sum of $9+7+5+5+1$ dimensional $\algg_{[0]}(={\mathfrak{sl}}(2,\C))$-modules.

Our next task is to compute leading terms of coefficients of twisted vertex operators associated to the elements of $R,$ or in
the other words, we have to describe the set $\leadt(\bar{R}^\sigma).$

\begin{lemma}\label{ltLemma} The set $\leadt(\bar{R}^\sigma)$ consists of the elements of the form:
\begin{eqnarray*}
& &X_{i_1}(j-\frac{1}{2})X_{i_2}(j-\frac{1}{2}),\ X_{i_1}(j-1)X_{i_2}(j)\\ & &\text{ with colors
}i_1i_2:11,21,22,31,32,41,33,42,51,43,52,44,53,54,55\\ & &X_{i_1}(j-\frac{3}{2})X_{i_2}(j-\frac{1}{2}),\
X_{i_1}(j-1)X_{i_2}(j-1)\\ & &\text{ with colors }i_1i_2:11,12,13,66,14,76,15,77,86,25,87,35,88,45,55\\ &
&X_{i_1}(j-\frac{1}{2})X_{i_2}(j)\\ & &\text{ with colors }i_1i_2:16,17,26,18,27,28,37,47,38,48,57,58\\ &
&X_{i_1}(j-1)X_{i_2}(j-\frac{1}{2})\\ & &\text{ with colors }i_1i_2:61,71,62,63,72,73,64,74,65,84,75,85\\
\end{eqnarray*}
where $j\in\Z.$
\end{lemma}
\begin{proof} Let us take the element $X_1(-1)X_6(-1)\vac,$ and let us demonstrate how one can compute
the leading terms for the corresponding vertex operator. 
Using the expression for the coefficient $(a(-1)b)(n)$ in the $(-1)$-product
(\cite{Li2}, definition 3.7) for $a=X_1(-1)\vac$ and $b=X_6(-1)\vac,$ and taking into account that $T=2,$ we
get
\begin{multline*}
Y^\sigma(X_1(-1)X_6(-1)\vac,z)=\sum_{i<0}X_1(i+\frac{1}{2})Y^\sigma(X_6(-1)\vac,z)z^{-i-1-\frac{1}{2}}+\\
+\sum_{i\geq 0}Y^\sigma(X_6(-1)\vac,z)X_1(i+\frac{1}{2})z^{-i-1-\frac{1}{2}}+\\ \sum_{k\geq
0}(-1)^k\binom{k-\frac{1}{2}}{k+1}Y^\sigma(X_1(k)X_6(-1)\vac,z)z^{-k-1}.
\end{multline*}
Since $X_1(k)X_6(-1)\vac=0$ for $k\geq 2$ on $N(\Lambda_0),$ we get
\begin{multline*}
Y^\sigma(X_1(-1)X_6(-1)\vac,z)=\sum_{i<0}\sum_{k\in\Z}X_1(i+\frac{1}{2})X_6(k)z^{-k-i-2-\frac{1}{2}}+\\
+\sum_{i\geq 0}\sum_{k\in\Z}X_6(k)X_1(i+\frac{1}{2})z^{-k-i-2-\frac{1}{2}}+\\
Y^\sigma([X_1,X_6](-1)\vac,z)z^{-1}-\frac{1}{2}\langle X_1,X_6\rangle
\end{multline*}
If we fix degree of $z$ to be $2j-\frac{1}{2}=-k-i-2-\frac{1}{2},$ we find that corresponding coefficient of
this vertex operator is of the form
\begin{multline*}
\sum_{i<0}X_1(i+\frac{1}{2})X_6(-2j-i-2)+\\ \sum_{i\geq 0}X_6(-2j-i-2)X_1(i+\frac{1}{2}) + \text{ ``larger
terms''}
\end{multline*}
where by ``larger terms'' we mean larger colored partitions. The definition of ordering and leading terms gives
us that the smallest colored partition in this sum is $X_1(j-\frac{1}{2})X_6(j).$

Repeating this process for all of 27 basis vectors of $R,$ and repeatedly taking into account leading terms we
have already calculated, we get the lemma.
\end{proof}

\subsection{Relations among relations}

Next, we count all possible embeddings of leading terms from Lemma \ref{ltLemma} in colored partitions of
length 3:

\begin{lemma}\label{EmbeddLemma} For $s\in\Z$ we have
\begin{eqnarray*}
& &\sum_{|\pi|=3s-\frac{5}{2}}N(\pi)=81,\ \sum_{|\pi|=3s-\frac{3}{2}}N(\pi)=81,\
\sum_{|\pi|=3s-\frac{1}{2}}N(\pi)=81\\ & &\sum_{|\pi|=3s-2}N(\pi)=80\ \sum_{|\pi|=3s-1}N(\pi)=80\
\sum_{|\pi|=3s}N(\pi)=80\\
\end{eqnarray*}
\end{lemma}
For example, colored partition $\pi=X_3(-\frac{3}{2})X_2(-\frac{3}{2})X_6(-1)$ has 2 embeddings
$\rho\subset\pi$: those are
$\rho=X_3(-\frac{3}{2})X_2(-\frac{3}{2})$ and $\rho=X_2(-\frac{3}{2})X_6(-1).$ On the other hand, colored partition
$X_3(-\frac{3}{2})X_6(-1)$ is not in $\leadt(\bar{R}^\sigma).$

Recall now the definition (\ref{defPhi}) of the map $\Phi.$ We have

\begin{lemma}\label{MPlemma} Let $\lbrace e_\theta=[e_1,e_2],e_1,e_2,h_1,h_2,f_1,f_2,[f_2,f_1]\rbrace$ be a
basis of $\algg$ build from the set of its canonical generators. Let $u_i,u^i, i=1,\ldots\dim\algg$ be a pair
of dual bases of $\algg$ with respect to $\langle\cdot,\cdot\rangle.$ The following vectors are highest weight
vectors of $\algg$-submodules of $N(\Lambda_0)\otimes N(\Lambda_0)$ of dimensions 64, 35, 35 and 27
respectively:
\begin{eqnarray*}
q_{64}&=&e_\theta(-2)\otimes e_\theta(-1)e_\theta(-1)-e_\theta(-1)\otimes e_\theta(-2)e_\theta(-1)\\
q_{35}&=&f_1(1)\cdot q_{64}=e_2(-1)\otimes e_\theta(-1)e_\theta(-1)-e_\theta(-1)\otimes e_2(-1)e_\theta(-1)\\
q_{\underline {35}}&=&-f_2(1)\cdot q_{64}= e_1(-1)\otimes e_\theta(-1)e_\theta(-1)-e_\theta(-1)\otimes
e_1(-1)e_\theta(-1)\\ q_{27}&=&\frac{1}{8}\sum_{i=1}^{\dim\algg}u_i(-1)\otimes u^i(0)e_\theta(-1)e_\theta(-1)-
2\cdot\vac\otimes e_\theta(-2)e_\theta(-1)
\end{eqnarray*}
Moreover, the corresponding $\algg$-submodules $Q_{64},Q_{35}, Q_{\underline{35}},Q_{27}$ are in $\ker\Phi$
(here we write $e_\theta(-1)e_\theta(-1)$ for $e_\theta(-1)e_\theta(-1)\vac$ etc. for short).
\end{lemma}
\begin{proof} It is straightforward to prove that vectors $q_{64},q_{35}$ and $q_{\underline{35}}$
are highest weight vectors for the action of $\algg,$ and that they lie in $\ker\Phi.$ Dimensions of the
corresponding modules are calculated using Weyl formula, and since $\Phi$ is a $\algg$-module map, it is also
clear that $Q_{64},Q_{35}$ and $Q_{\underline{35}}$ are in $\ker\Phi.$

Let $x\in\algg.$ Write
$$[x,u_i]=\sum_j\alpha_{ij}u_j,\ [x,u^i]=\sum_j\beta_{ji}u^j.$$
By taking the product of the first equation with $\langle\cdot,u^j\rangle$ and second with
$\langle\cdot,u_j\rangle$ we get
$$\alpha_{ij}=\langle[x,u_i],u^j\rangle,\ \beta_{ji}=\langle[x,u^i],u_j\rangle.$$
Using invariance of $\langle\cdot,\cdot\rangle$ we obtain
$$\alpha_{ij}=\langle x,[u_i,u^j]\rangle,\ \beta_{ji}=\langle x,[u^i,u_j]\rangle$$
so it follows $\alpha_{ij}=-\beta_{ji}.$ Now we get
\begin{eqnarray*}
[x(0),q_{27}]&=&\frac{1}{8}\sum_{i=1}^{\dim\algg}[x(0),u_i(-1)]\otimes u^i(0)e_\theta(-1)e_\theta(-1)+\\ &
&\frac{1}{8}\sum_{i=1}^{\dim\algg}u_i(-1)\otimes [x(0),u^i(0)]e_\theta(-1)e_\theta(-1)+\\ &
&\frac{1}{8}\sum_{i=1}^{\dim\algg}u_i(-1)\otimes u^i(0)x(0)e_\theta(-1)e_\theta(-1)+\\ & &2\cdot\vac\otimes
x(0)e_\theta(-2)e_\theta(-1)
\end{eqnarray*}
By using $\alpha_{ij}=-\beta_{ji}$ it follows that
\begin{eqnarray*}
[x(0),q_{27}]&=&\frac{1}{8}\sum_{i=1}^{\dim\algg}u_i(-1)\otimes u^i(0)x(0)e_\theta(-1)e_\theta(-1)+\\ &
&2\cdot\vac\otimes x(0)e_\theta(-2)e_\theta(-1)
\end{eqnarray*}
Since $e_\theta(-1)e_\theta(-1)$ and $e_\theta(-2)e_\theta(-1)$ are highest weight vectors, we get
$$[e_1(0),q_{27}]=[e_2(0),q_{27}]=[e_\theta(0),q_{27}]=0$$
as required.

Recall now that Virasoro element $\omega$ of $N(\Lambda_0)$ is defined by
$$\omega=\frac{1}{8}\sum_{i=1}^{\dim\algg}u_i(-1)u^i(-1).$$
Derivation of the vertex operator algebra $N(\Lambda_0)$ is then
$$L(-1)=\omega(0)=\frac{1}{8}\sum_{i=1}^{\dim\algg}\sum_{n+m=-1}\NO u_i(n)u^i(m)\NO.$$
Since $\algg(n)e_\theta(-1)e_\theta(-1)=0$ for $n>0,$ we get
\begin{eqnarray*}
\Phi(q_{27})&=&\frac{1}{8}\sum_{i=1}^{\dim\algg}u_i(-1)u^i(0)e_\theta(-1)e_\theta(-1)-
e_\theta(-2)e_\theta(-1)=\\ &=&L(-1)e_\theta(-1)e_\theta(-1)-2e_\theta(-2)e_\theta(-1)=0
\end{eqnarray*}
\end{proof}

In order to apply Theorem \ref{theThm}, we need to determine $\dim \sum Q_p(n),p=64,35,\underline{35},27.$ 
We shall use the irreducibility of loop $\tilde{\algg}[\sigma]$-modules
arising from irreducible $(\algg,\sigma)$-modules established in Proposition \ref{realLoopIrrdb}.

\begin{lemma}$\algg$-submodules $Q_{64},Q_{35}, Q_{\underline{35}},Q_{27}$ are irreducible
$(\algg,\sigma)$-modules.
\end{lemma}
\begin{proof}First, since $Q_{64}$ is an irreducible $\algg$-module and $\sigma(q_{64})=q_{64},$ it
follows that $Q_{64}$ is an irreducible $(\algg,\sigma)$-module. Second, because
$\sigma(q_{35})=q_{\underline{35}},$ and $Q_{35}$ and $Q_{\underline{35}}$ are irreducible $\algg$-modules,
$Q_{35}\oplus Q_{\underline{35}}$ is also an irreducible $(\algg,\sigma)$-module. Finally, if $(u_i,u^i)$ is
pair of elements of dual bases, then they lie in the same $\sigma$-eigenspace of $\algg,$ so
$\sigma(q_{27})=q_{27},$ providing that $Q_{27}$ is closed under $\sigma,$ and therefore $Q_{27}$ is also an
irreducible $(\algg,\sigma)$-module.
\end{proof}
Let us temporarily denote by $Q\subset N(k\Lambda_0)\otimes N(k\Lambda_0)$ {\em any} $(\algg,\sigma)$-module
invariant under $L(0)$ as well. Introduce the notation
$$\tilde{\algg}_{\geq 0}=\bigoplus_{n\geq 0}\algg(n),\ \tilde{\algg}_{> 0}=\bigoplus_{n> 0}\algg(n)$$
$$Q_{\geq 0}=U(\tilde{\algg}_{\geq 0})Q,\ Q_{> 0}=U(\tilde{\algg}_{> 0})Q$$

Since $Q$ is a $(\algg,\sigma)$-module, we may construct the twisted loop module for $\tilde{\algg}[\sigma]$ from $Q$
as in section \ref{loopSection}. Denote that module by $Q^\sigma_A.$

Set (cf. (\ref{tensorCoef}))
$$\overline{(Q_{\geq 0})^\sigma}=\coprod_{n\in\frac{1}{T}\Z}Q_{\geq 0}(n)\subset U^\sigma\bar{\otimes}U^\sigma$$
$$\overline{(Q_{>0})^\sigma}=\coprod_{n\in\frac{1}{T}\Z}Q_{>0}(n)\subset U^\sigma\bar{\otimes}U^\sigma$$

\begin{proposition}\label{loopIso} The map
$$\phi:Q^\sigma_A\rightarrow
\frac{{\overline{(Q_{\geq 0})^\sigma}}}{{\overline{(Q_{> 0})^\sigma}}} ,\ \ \phi(q\otimes t^n)=
q(n)+Q_{>0}(n)$$ is homomorphism of $\frac{1}{T}\Z$-graded $\tilde{\algg}[\sigma]$-modules.
\end{proposition}
\begin{proof} The twisted vertex operator commutator formula gives
\begin{equation}\label{twBracket}
[x(m),q(n)]=\sum_{i\geq 0}\binom{m}{i}(x(i)q)(m+n-i)
\end{equation}
for all $x\in\algg,q\in Q$ and all $m,n\in\frac{1}{T}\Z.$ Since $Q_{>0}$ is $\tilde{\algg}_{\geq 0}$-module, it
follows from (\ref{twBracket}) that
\begin{equation}\label{loopOperDef}
[x(m),q(n)+Q_{>0}(n)]=(x(0)q)(m+n)+Q_{>0}(m+n)
\end{equation}
providing that the quotient
$$\frac{{\overline{(Q_{\geq 0})^\sigma}}}{{\overline{(Q_{> 0})^\sigma}}}$$
is well-defined $\tilde{\algg}[\sigma]$-module. Comparing (\ref{loopAbstrDef}) with (\ref{loopOperDef})
together with identification $L(0)\leftrightarrow -d$ shows that $\phi$ is indeed homomorphism of
$\tilde{\algg}[\sigma]$-modules, so the proof is complete.
\end{proof}
\begin{remark}\label{remPhiIn} By using Propositions \ref{realLoopIrrdb} and \ref{loopIso},
we see that if $Q$ is $(\algg,\sigma)$-irreducible, then $\phi$ is either zero or it is injective.
\end{remark}
\begin{remark}\label{remPhiDim} Let $n\in\Z$ be a fixed integer, and let $Q\subset N(k\Lambda_0)\otimes N(k\Lambda_0)$
be a finite dimensional $(\algg, \sigma)$-irreducible $L(0)$-invariant module. Then
$$\dim(Q)=\sum_{i=0}^{T-1}\dim\left(Q^\sigma_A(n+\frac{i}{T})\right).$$
Suppose that $\phi$ is not zero. Then Proposition \ref{loopIso} and Remark \ref{remPhiIn} imply that
$$\dim(Q)=\sum_{i=0}^{T-1}\dim\left(Q(n+\frac{i}{T})\right).$$
In the other words, map $q\mapsto q(n)$ is injective on $Q.$
\end{remark}

\begin{lemma} Let $Q_{64},Q_{35}, Q_{\underline{35}},Q_{27}$ be as in the previous lemma. Then
for any $n\in\Z$ and $p=64,35,\underline{35},27$
\begin{eqnarray*}
\dim Q_{p}(n)&=&\sum_{i=0,1}\dim\left(Q(n+\frac{i}{T})\right)=\dim Q_{p}
\end{eqnarray*}
\end{lemma}
\begin{proof} By remark \ref{remPhiDim}, we have to prove that map $\phi$ is not identically zero. An easy calculation
shows that
\begin{eqnarray*}
(Q_{64})_{>0}&=&Q_{35}\oplus Q_{\underline{35}}\oplus \left(\univ(\algg)\cdot(\vac\otimes
e_\theta(-1)e_\theta(-1))\right)
\end{eqnarray*}
Denote by ${E_1,H_1,F_1}$ (cf. section 3.1) the basis of Lie algebra $\algg_{[0]},$ and let $\alpha$ be its
simple root.

For $q_{64}$ and for all $n\in\Z$ we have
$$H_1(0)\cdot q_{64}(n)=6\alpha(H_1)q_{64}(n).$$
On the other hand, the vector in $(Q_{64})_{>0}(n)$ with highest possible weight for the action of
$\algg_{[0]}$ is $q_{35}+q_{\underline{35}},$ and its weight is
$$H_1(0)\cdot (q_{35}+q_{\underline{35}})(n)=5\alpha(H_1)(q_{35}+q_{\underline{35}})(n).$$
This calculation shows that $q_{64}(n)+(Q_{64})_{>0}(n)\neq (Q_{64})_{>0}(n),$ i.e. $\phi$ is nonzero. The same
argument applies to the remaining cases, so the proof is complete.
\end{proof}
Since $(Q_{35}+Q_{\underline{35}})(n)\subseteq ((Q_{64})_{>0})(n)$ it follows that 
$Q_{64}(n)\cap (Q_{35}+Q_{\underline{35}})(n)=0$. Likewise, $Q_{27}(n)\subseteq ((Q_{35}+Q_{\underline{35}})_{>0})(n),$
so the sum 
$$ Q_{64}(n)\oplus (Q_{35}+Q_{\underline{35}})(n)\oplus Q_{27}(n)$$
is indeed direct. Therefore, for
$$Q=Q_{64}\oplus Q_{35}\oplus Q_{\underline{35}}\oplus Q_{27}$$
and arbitrary $n\in\Z$ holds that
\begin{eqnarray}
\dim Q(n)=161\label{myHRel}
\end{eqnarray}
\subsection{Initial conditions}

Let $\Lambda_0,\Lambda_1$ be the fundamental weights associated to the chosen Cartan subalgebra of
$\tilde{\algg}[\sigma].$ Basic module is the (unique) standard module with the highest weight
$\Lambda=0\cdot\Lambda_0+1\cdot\Lambda_1$. Denote by $v_\Lambda$ highest weight vector of basic module. The
following relations on the basic module hold:
$$e_0\cdot v_\Lambda=e_1\cdot v_\Lambda=0,\ f_0\cdot v_\Lambda=0,{f_1}^2\cdot v_\Lambda=0,\ 
h_0\cdot v_\Lambda=0,\ h_1\cdot v_\Lambda=v_\Lambda.$$ In our notation it means that
\begin{eqnarray*}
& &X_5(\frac{1}{2})\cdot v_\Lambda= X_6(0)\cdot v_\Lambda=0\\ & &X_1(-\frac{1}{2})\cdot v_\Lambda =0,\
X_8(0)X_8(0)\cdot v_\Lambda=0\\ & &X_7(0)\cdot v_\Lambda=v_\Lambda
\end{eqnarray*}

\subsection{The main result}

Recall that we denoted by ${\cal P}$ a basis for ${\cal S}^\sigma=S(\bar{\algg}[\sigma]).$ We shall introduce the
following notation in our special case $\tilde{\algg}[\sigma]=A_2^{(2)}$:
$$S^\sigma_-=S(\bar{\algg}[\sigma]_{<0}+\C X_8(0)+\C X_1(-\frac{1}{2})).$$
Let ${\cal P}_-$ be a basis for $S^\sigma_-.$ 

We shall say that partition 
$\pi=X_{i_1}(j_1)X_{i_2}(j_2)\cdots X_{i_n}(j_n)\in {\cal P}_-$ satisfies the difference
${\cal R}$ condition if it does not contain any of partitions from the set $\leadt(\bar{R}_\sigma).$
\begin{theorem}\label{level1Basis} The set of vectors
\begin{equation}\label{reducedSpan}
u(X_{i_1}(j_1)X_{i_2}(j_2)\cdots X_{i_n}(j_n))v_\Lambda,
\end{equation}
where $\pi=X_{i_1}(j_1)X_{i_2}(j_2)\cdots X_{i_n}(j_n)\in{\cal P}_-$
satisfies the difference $\cal{R}$ condition and {\em initial} condition $X_1(-\frac{1}{2})\not\subset \pi$
is a basis of the basic $A_2^{(2)}$-module $L(\Lambda_1)$.
\end{theorem}
\begin{proof} By using lemma \ref{EmbeddLemma} and (\ref{myHRel}) we see that the conditions of the Theorem
\ref{theThm} are fulfilled. Denote by $\rho(\pi)$ exactly one of all possible
embeddings $\rho\subset\pi,\ \rho\in\leadt(\bar{R}_\sigma),\ \pi\in{\cal P}_-.$ 
Then set
$$u(\rho(\pi)\subset \pi)v_{\Lambda_1},\ \pi\in{\cal P}_-$$
is a basis for the maximal submodule $M^1(\Lambda_1)$ 
if {\em initial} conditions determined by module definition relations 
$h_0\cdot v_\Lambda=0,\ h_1\cdot v_\Lambda=v_\Lambda$ produce new leading terms for which
the statement of the Theorem \ref{theThm} still holds. For example, from 
$X_1(-\frac{1}{2})v_{\Lambda_1}=0$ we get that $\rho'=X_1(-\frac{1}{2})$ is new leading term. Now we have to check
that for any $\rho\in\leadt(\bar{R}^\sigma)$ and $\pi\in{\cal P}_-$ such that $\rho',\rho\subset\pi$ there is a
relation among relations between $u(\rho'\subset\pi)$ and $u(\rho\subset\pi).$

Calculation by hand shows that this is indeed true for any new leading term produced by module definition relations,
so finally by using procedure described in section 1.7 we obtain that
(\ref{reducedSpan}) is a basis for $L(\Lambda_1).$
\end{proof}

Reformulating description of $\leadt(\bar{R}_\sigma)$, we may say that if $u(\pi)v_{\Lambda_1}$ from (\ref{reducedSpan}) is
in the basis of $L(\Lambda_1),$ then any $\rho\subset\pi,$ whose shape is described below, must also satisfy
colors condition of the form:
\begin{eqnarray*}
& &X_{i_1}(j-\frac{1}{2})X_{i_2}(j-\frac{1}{2}),\ X_{i_1}(j-1)X_{i_2}(j),\ i_1,i_2\in\lbrace 2,3,4\rbrace,\
i_2\geq i_1\\ & &X_{i_1}(j-\frac{3}{2})X_{i_2}(j-\frac{1}{2}),\ X_{i_1}(j-1)X_{i_2}(j-1),\ i_1,i_2\in\lbrace
6,7,8\rbrace,\ i_2\geq i_1\\ & &X_{i_1}(j-\frac{1}{2})X_{i_2}(j),\ i_1\in\lbrace 3,4,5\rbrace,\ i_2=8\\ &
&X_{i_1}(j-1)X_{i_2}(j-\frac{1}{2}),\ i_1=8,i_2\in\lbrace 1,2,3\rbrace\\
\end{eqnarray*}
where $j\in\Z.$

\section{Combinatorial identities}
Let $(s_0,s_1)$ be a pair of nonnegative relatively prime integers. The specialization of type ${\bf s}=(s_0,s_1)$
amounts to the following:
\begin{eqnarray*}
& &X_1(-k-\frac{1}{2})\mapsto s_0+k(4s_1+2s_0)\\ & &X_2(-k-\frac{1}{2})\mapsto s_0+s_1+k(4s_1+2s_0)\\ &
&X_3(-k-\frac{1}{2})\mapsto s_0+2s_1+k(4s_1+2s_0)\\ & &X_4(-k-\frac{1}{2})\mapsto s_0+3s_1+k(4s_1+2s_0)\\ &
&X_5(-k-\frac{1}{2})\mapsto s_0+4s_1+k(4s_1+2s_0)\\ & &X_6(-k-1)\mapsto -s_1+(k+1)(4s_1+2s_0)\\ &
&X_7(-k-1)\mapsto (k+1)(4s_1+2s_0)\\ & &X_8(-k)\mapsto s_1+k(4s_1+2s_0)
\end{eqnarray*}
for $k\geq 0.$
\begin{theorem}\label{thmWithColors} The number of partitions $(m_1^{(c_1)},\ldots,m_r^{(c_r)})$ of an integer $n$ in two
colors (i.e. $c_i=1,2$ for $i=1,\ldots,r$) such that difference between two consecutive parts is at least 4
(i.e. $m_{i+1}^{(c_{i+1})}-m_i^{(c_i)}\geq 4$) and
\begin{eqnarray*}
& &m_{i}^{(2)}\equiv\pm 1\mod 6,\ m_i^{(2)}> 1\text { for } i=1,\ldots,r\\ &
&m_{i+1}^{(c_{i+1})}-m_i^{(c_i)}=6\Rightarrow c_{i+1}\neq 2\text{ and } c_i\neq 2\\ &
&m_{i+1}^{(c_{i+1})}-m_i^{(c_i)}=5\Rightarrow c_{i+1}\neq 2\text{ or } c_i\neq 2\\ &
&m_{i+1}^{(c_{i+1})}-m_i^{(c_i)}=4\Rightarrow c_{i+1}\neq 2\text{ or } c_i\neq 2\text{ and }\\ & &\hskip 3.7cm
m_{i+1}^{(c_{i+1})}+m_i^{(c_i)}\not\equiv \pm 4\!\!\mod 12
\end{eqnarray*}
is equal to the number of partitions of $n$ in one color into parts congruent to $\pm 1$ modulo $6.$
\end{theorem}
\begin{proof} The specialization of type ${\bf s}=(1,1)$ gives the character formula
$$\ch_{(1,1)}L(\Lambda_1)=\prod_{n\equiv\pm 1\!\!\!\mod 6}\frac{1}{1-q^n}.$$
This specialization also gives
\begin{eqnarray*}
& &X_1(-k-\frac{1}{2})\mapsto (6k+1)^{(2)},\hskip1cm X_6(-k-1)\mapsto (6k+5)^{(1)}\\ &
&X_2(-k-\frac{1}{2})\mapsto (6k+2)^{(1)},\hskip1cm X_7(-k-1)\mapsto (6k)^{(1)}\\ & &X_3(-k-\frac{1}{2})\mapsto
(6k+3)^{(1)},\hskip1cm X_8(-k)\mapsto (6k+1)^{(1)}\\ & &X_4(-k-\frac{1}{2})\mapsto (6k+4)^{(1)}\\ &
&X_5(-k-\frac{1}{2})\mapsto (6k+5)^{(2)}
\end{eqnarray*}
for $k\geq 0.$ If we interpret difference ${\cal R}$ condition and initial condition 
in terms of this specialization, we get the list of forbidden parts in partitions. By forming a table with allowed
differences and colors between consecutive parts, we get the Theorem.
\end{proof}
{\bf Proof of Theorem A:}
The specialization of type ${\bf s}=(2,1)$ gives the character formula
$$\ch_{(2,1)}L(\Lambda_1)=\prod_{n\geq 1}(1+q^{2n-1}).$$
On the other hand, we have that
\begin{eqnarray*}
& &X_1(-k-\frac{1}{2})\mapsto 8k+2,\hskip1cm X_6(-k-1)\mapsto 8k+7\\ & &X_2(-k-\frac{1}{2})\mapsto
8k+3,\hskip1cm X_7(-k-1)\mapsto 8k\\ & &X_3(-k-\frac{1}{2})\mapsto 8k+4,\hskip1cm X_8(-k)\mapsto 8k+1\\ &
&X_4(-k-\frac{1}{2})\mapsto 8k+5\\ & &X_5(-k-\frac{1}{2})\mapsto 8k+6
\end{eqnarray*}
Again, interpreting difference ${\cal R}$ condition and initial condition in the same way we did it in the proof of 
Theorem \ref{thmWithColors}, we get the statement.

\vskip 1cm
{\small\em
\noindent
Ivica Siladi\'c\\
Department of Mathematics\\
University of Zagreb\\
Bijeni\v{c}ka 30, 10000 Zagreb\\
Croatia\\
e-mail: siladic@math.hr
}
\end{document}